\documentclass{article}
\usepackage{amsmath, epsfig,amssymb, multicol}
\usepackage{setspace}
\newtheorem{lemma}{Lemma}
\newtheorem{theorem}{Theorem}

\newtheorem{question}{Qeustion}

\newtheorem{fact}{Fact}

\newcommand{\fr}{{\cal F}_{\rm{red}}}
\newcommand{\fb}{{\cal F}_{\rm{blue}}}
\newcommand{\pr}{{\bf Property}}
\newcommand{\p}{{\cal P}}
\renewcommand{\c}{{\cal C}}

\title{The Ramsey number of generalized loose paths in  hypergraphs}
\author{Xing Peng
\thanks{Center for Applied Mathematics, Tianjin University, 300072, P.~R.~China,
({\tt pengxingmath@gmail.com}).}}
\date{}
\begin{document}
\maketitle

\begin{abstract}
Let $H=(V,E)$ be an $r$-uniform hypergraph.  For each $1 \leq s \leq r-1$, an $s$-path $\p^{r,s}_n$   of length $n$ in $H$  is a sequence of distinct vertices $v_1,v_2,\ldots,v_{s+n(r-s)}$ such that $\{v_{1+i(r-s)},\ldots, v_{s+(i+1)(r-s)}\}\in E(H)$  for each $0 \leq i \leq n-1$.
Recently, the Ramsey number of  $1$-paths  in uniform hypergraphs has received a lot of attention.
 In this paper, we   consider the Ramsey number of  $r/2-$paths for even $r$.  Namely, we prove the following exact result:
$R(\p^{r,r/2}_n,\p^{r,r/2}_3)=R(\p^{r,r/2}_n,\p^{r,r/2}_4)=\tfrac{(n+1)r}{2}+1.$
\end{abstract}

\section{Introduction}
An $r$-uniform hypergraph $H$ is a pair $H=(V,E)$, where $V$ is a set of vertices and $E$ is a collection of $r$-subsets of $V$.
For two $r$-uniform hypergraphs $H_1$ and $H_2$,  the {\it Ramsey number} $R(H_1,H_2)$ is the minimum value of $N$ such that each red-blue coloring of edges in the complete $r$-uniform hypergraph $K_N^r$ on $N$ vertices  contains either  a red $H_1$  or  a blue $H_2$.
 Let $H$ be  an $r$-uniform hypergraph. For each $1 \leq s \leq r-1$, an {\it $s$-path} $\p^{r,s}_n$ of length $n$ in $H$  is a sequence of distinct vertices $v_1,v_2,\ldots,v_{s+n(r-s)}$ such that
 $\{v_{1+i(r-s)},\ldots, v_{s+(i+1)(r-s)}\}$ is an edge of $H$ for each $0 \leq i \leq n-1$. Similarly, an {\it $s$-cycle} $\c^{r,s}_n$ of length $n$ in $H$ is a sequence of vertices $v_1,v_2,\ldots,v_{s+n(r-s)}$ such that $\{v_{1+i(r-s)},\ldots, v_{s+(i+1)(r-s)} \}$ is an edge of $H$ for each $0 \leq i \leq n-1$, $v_1,\ldots,v_{n(r-s)}$ are distinct, and $v_{n(r-s)+j}=v_j$ for each $1 \leq j \leq s$. An $s$-path (resp. an $s$-cycle) is {\it loose} if $1 \leq s \leq r/2$ and  an  $s$-path  (resp. an $s$-cycle) is {\it tight} if $r/2 < s \leq r-1$.

When $r=2$ and $s=1$, we get  paths and cycles in graphs. A classical result from Ramsey theory \cite{gg} says $R(P_n,P_m)=n+\lfloor \tfrac{m+1}{2}\rfloor$ for $n \geq m \geq 1$; it is also known \cite{flp,fs} that $R(P_n,C_m)=R(P_n,P_m)=n+\tfrac{m}{2}$ for $n \geq m$ and $m$ even. One may ask what is the Ramsey number of paths and cycles in uniform hypergraphs?

The following construction \cite{ha} was used   to show a lower bound on $R(\p^{3,1}_n, \p^{3,1}_n)$ for $n \geq 1$; we can adapt it to show   $R(\p^{r,s}_n, \p^{r,s}_m) > s+n(r-s)+\lfloor \tfrac{m+1}{2} \rfloor-2$ for $n \geq m \geq 1$ and $1 \leq s \leq r/2$.
 To see this, we let $N=s+n(r-s)+\lfloor \tfrac{m+1}{2} \rfloor-2$ and partition the vertex set of $K^r_N$ into two subsets $A$ and $B$, where $|A|=s+n(r-s)-1$ and $|B|=\lfloor \tfrac{m+1}{2} \rfloor-1$.  We color all edges $f$  satisfying  $f \subseteq A$ or $f \subseteq B$ red and the remaining edges blue. Observe that the number of vertices in an $s$-path with length $n$ equals $s+n(r-s)$, so there is no red $\p^{r,s}_n$. Since each vertex in a loose path can be in at most two  edges, a blue path $\p^{r,s}_m$ must have at least $\lfloor \frac{m+1}{2} \rfloor$ vertices from $B$. As the assumption on $|B|$, there is no blue  $\p^{r,s}_m$. We showed that $N=s+n(r-s)+\lfloor \tfrac{m+1}{2}\rfloor-2$ is a lower bound for $R(\p_n^{r,s},\p_m^{r,s})$.

We have the following  interesting  question which asks whether the construction above gives the true value of $R(\p^{r,s}_n,\p^{r,s}_m)$.
 \begin{question}\label{q:q1}
 Is $R(\p^{r,s}_n, \p^{r,s}_m)=s+n(r-s)+\lfloor \tfrac{m+1}{2} \rfloor-1$ for $n \geq m \geq 1$ and $1 \leq s \leq r/2$?
 \end{question}
 This question can be answered in the affirmative way for   the case where $s=1$.  Haxell et al. \cite{ha} first determined the asymptotic values of  $R(\p^{3,1}_n,\p^{3,1}_n)$,  $R(\c^{3,1}_n,\c^{3,1}_n)$, and $R(\p^{3,1}_n,\c^{3,1}_n)$.  Later,
   Gy\'arf\'as, S\'ark\"ozy, and Szemer\'edi \cite{grs}
   extended this result to all $r \geq 3$. Namely, they proved that $R(\p^{r,1}_n,\p^{r,1}_n)$,  $R(\p^{r,1}_n,\c^{r,1}_n)$, and $R(\c^{r,1}_n,\c^{r,1}_n)$ are asymptotically equal to $\tfrac{(2r-1)n}{2}$.  There are some exact results on short paths and cycles.  Gy\'arf\'as and Raeisi \cite{gr} proved
   \[ R(\p^{r,1}_3,\p^{r,1}_3)=R(\p^{r,1}_3,\c^{r,1}_3)=R(\c^{r,1}_3,\c^{r,1}_3)+1=3r-1;\]
in the same paper,  they also proved
  \[R(\p^{r,1}_4,\p^{r,1}_4)=R(\p^{r,1}_4,\c^{r,1}_4)=R(\c^{r,1}_4,\c^{r,1}_4)+1=4r-2.\]
For $r=3$ and $s=1$,  Maherani et al.\cite{m1} determined the exact value of $R(\p^{3,1}_n,\p^{3,1}_m)$
for $n \geq \lfloor \tfrac{5m}{4}\rfloor$.
Recently, Omidi and Shahsiah \cite{m2}  proved the following general result.  For $n \geq m \geq 1$, we have
\[
R(\p^{3,1}_n,\p^{3,1}_m)=R(\p^{3,1}_n,\c^{3,1}_m)=R(\c^{3,1}_n,\c^{3,1}_m)+1=2n+\lfloor \tfrac{m+1}{2} \rfloor \textrm{ and } R(\p^{3,1}_m,\c^{3,1}_n)=2n+\lfloor \tfrac{m-1}{2} \rfloor.
\]
For more details on small Ramsey numbers, the reader is referred to the dynamic survey paper \cite{r1}.

  To the author's best knowledge, there is no attempt to study the Ramsey number of other types of paths in hypergraphs.  In this paper, we will show some exact results for $s=r/2$ and $r$ even. As we will see, the following lemma will be important for establishing these results.
\begin{lemma}
For each $s \geq 1$ and $n \geq 2$, we have
 \[R(\p_n^{2s,s},\p_2^{2s,s})=(n+1)s.\]
\end{lemma}

\noindent
{\bf Proof:} We will prove the lemma by induction on $n$. It is easy to see that
 $R(\mathcal{P}^{2s,s}_{2},\mathcal{P}^{2s,s}_{2})=3s$. Now let $n\geq
 3$ and $c$ be a red-blue coloring of the edges in $K_N^{2s}$, where $N=(n+1)s$. By  induction hypothesis we have
 $R(\mathcal{P}^{2s,s}_{n-1},\mathcal{P}^{2s,s}_{2})<(n+1)s$.
 So there is either a red $\mathcal{P}^{2s,s}_{n-1}$ or a blue
 $\mathcal{P}^{2s,s}_{2}$. We need only consider the former case.
 Assume $A_1,A_2,\ldots,A_{n}$ is a red
 $\mathcal{P}^{2s,s}_{n-1}$ and $A_{n+1}$ is the remaining $s$ vertices. Consider the edges $g=A_1\cup A_{n+1}$
 and $h=A_n\cup A_{n+1}$. If at least one of $g$ or $h$ is red, then we
 have a red $\mathcal{P}^{2s,s}_{n}$. Otherwise, we have a blue
 $\mathcal{P}^{2s,s}_{2}$.
\hfill $\square$

 We will prove the following two main theorems.
\begin{theorem}\label{t:1}
For each $s \geq 1$ and $n \geq 3$, we have
 \[R(\p_n^{2s,s},\p_3^{2s,s})=(n+1)s+1.\]
\end{theorem}

\begin{theorem}\label{t:2}
For each $s \geq 1$ and $n \geq 4$, we have
  \[
  R(\p_n^{2s,s},\p_4^{2s,s})=(n+1)s+1.
  \]
\end{theorem}
Notice that theorems above provide a partial positive answer to Question \ref{q:q1} for $s=r/2$ and $r$ even. To prove Theorem \ref{t:1} and Theorem \ref{t:2}, we will need only prove the upper bound.

Throughout this paper, for a red-blue coloring of edges in a uniform hypergraph, we use $\fr$ (resp. $\fb$) to denote the subhypergraph induced by all red (resp.  blue) edges. For a positive integer $N$, we use $[N]$ to denote the set of the first $N$ positive integers.
 Since we will work on a fixed type of path $\p^{2s,s}_n$ in Section 2 and in Section 3, we will drop the superscripts and write $\p_n$ for $\p^{2s,s}_n$ in these two sections. For $m \geq 2$,   let $\{A_1,A_2,\ldots,A_m\}$ be a collection of pairwise disjoint $s$-sets.
 If  $A_i \cup A_{i+1} $ is an edge  $f_i$ for $1 \leq i \leq m-1$,  then we will say $A_1,A_2,\ldots,A_m$ form an $s$-path of length $m-1$. In this case, we will also say $f_1,f_2,\ldots,f_{m-1}$ induce this $s$-path of length $m-1$. We will refer $A_1$ and $A_m$ as the ending $s$-sets of this red path.


The paper is organized as follows. Since the proof of Theorem \ref{t:2} requires Theorem \ref{t:1},
  we will prove Theorem \ref{t:1} in Section 2 and Theorem \ref{t:2}  in Section 3. We will give some concluding remarks in the last section.

\section{Proof of Theorem \ref{t:1}}
For a fixed $s$, Theorem \ref{t:1} will be proved by induction on $n$. Because the idea for proving the base case and the inductive step are similar, we give an outline for the inductive step here. Suppose Theorem \ref{t:1} holds for all $3 \leq n \leq m-1$.  Let $c$ be a red-blue coloring of edges in $K_{(m+1)s+1}^{2s}$.
 By the inductive hypothesis, we have $(m+1)s+1 > R(\p_{m-1},\p_3)$. Thus either there is  a red $\p_{m-1}$, or there is a blue $\p_3$.  We need only consider the former case and also assume that there is no red $\p_m$. Let $\{A_1,A_2,\ldots,A_m\}$ be a family of  pairwise disjoint $s$-sets of $[ms]$ and $B$ be the remaining $s+1$ vertices in $[(m+1)s+1]$.
 We assume a red path is $A_1,A_2,\ldots,A_m$ and aim to find a blue $\p_3$. We write the edge with vertices $A_i \cup A_{i+1}$ as $f_i$ for each $1 \leq i \leq m-1$.
 For each $0 \leq l \leq s$, we say an edge $f$ is of type $(l,s-l,s)$ if $|f \cap B|=l$, $|f \cap A_2|=s-l$, and $A_p \subset f$ for some $1 \leq p  \leq m$ with $p \not =2$.
 We will couple  edges of type $(s+1-l,l-1,s)$ with edges of type $(l,s-l,s)$  as well as edges of type
 $(l,s-l,s)$ with edges of type  $(s-l,l,s)$. Lemma \ref{l:lm1} and Lemma \ref{l:lm2}  will show how  the color of edges of the first type forces the color of edges of the second type under some assumptions. We  note that $m$ is fixed and $m \geq 3$.

\begin{lemma} \label{l:lm1}
  Assume the edge $A_i \cup A_j $ is red for all $1 \leq i \not = j
  \leq m$ and there is no red $\p_m$.  For a fixed $1 \leq l \leq \lfloor \tfrac{s}{2} \rfloor$, if  all edges of type $(s+1-l,l-1,s)$ are blue and there  exists  a blue edge of   type $(l,s-l,s)$, then there is  a blue $\p_3$.
\end{lemma}
\noindent
{\bf Proof:} Assume that there is a blue edge $g_1$ of type $(l,s-l,s)$. We can assume $g_1=B'\cup A_2' \cup A_i$, where   $B'$ is an $l$-subset of $B$,  $A_2'$ is an $(s-l)$-subset of $A_2$, and $i \not =2$. Choose $A_2''$ to be an  $(l-1)$-subset of $A_2 \setminus A_2'$. Let $j \in \{1,m\} \setminus \{i\}$.
  We define
\[ g_2= (B \setminus B') \cup  A_2'' \cup A_j \textrm{ and } g_3= (B\setminus B') \cup  A_2'' \cup A_i.\]
We observe that both $g_2$ and $g_3$ are of type $(s+1-l,l-1,s)$.
By assumption, we get both $g_2$ and $g_3$ are blue.
Thus edges $g_1,g_3,$ and $g_2$ form a blue $\p_3$. \hfill $\square$

\begin{lemma}\label{l:lm2}
   Assume the edge $A_i \cup A_j$ is red   for all $1 \leq i \not = j
  \leq m$ and there is no red $\p_m$.  For a fixed $1 \leq l \leq \lfloor \tfrac{s}{2} \rfloor$, if  all edges of type $(l,s-l,s)$ are red, then  all edges of type $(s-l,l,s)$ are blue.
\end{lemma}
\noindent
{\bf Proof:} Suppose  that there is a red edge $g_1$ of type $(s-l,l,s)$. Without loss of generality, we can assume $g_1=B' \cup A_2' \cup A_i$, where   $B' \subseteq B $ with $|B'|=s-l$,  $A_2' \subseteq A_2$ with $|A_2'|=l$, and   $i \not =2$,  We pick an arbitrary $l$-subset $B''$ of $B \setminus B'$.

If $i=1$, then we  define
\[
 g_2= B'' \cup  (A_2 \setminus A_2') \cup A_1  \text{ and }  g_3=B''\cup  (A_2 \setminus A_2') \cup  A_3.
 \]
We notice  that both $g_2$ and $g_3$ are of type $(l,s-l,s)$ and they  are red  by assumption.   Now  $g_1,g_2,g_3,f_3,\ldots,f_{m-1}$ form a red $\p_m$, which is a contradiction .

If $i=3$, then we define
 \[
 g_2=B'' \cup (A_2 \setminus A_2') \cup A_3, \quad  g_3=B'' \cup (A_2 \setminus A_2') \cup  A_1, \text{ and } g_4=A_1 \cup A_4.
 \]
We get  $g_2$ and $g_3$ are red as they are of type $(l,s-l,s)$. The edge $g_4$ is also red by assumption,.  Now, $g_1,g_2,g_3,g_4,f_4,\ldots,f_{m-1}$ is a red $\p_m$, which is a contradiction.  We notice that if $m=3$, then we do not have to define $g_4$ as a red $\p_3$ with edges $g_1,g_2,g_3$ is sufficient.

If $i \not \in \{1,3\}$, then
  edges $B'' \cup (A_2 \setminus A_2') \cup A_1$ and $B'' \cup (A_2 \setminus A_2') \cup A_3$ are both red as they are of type $(l,s-l,s)$. Now $A_1,B''\cup (A_2\setminus A_2'),A_3,\ldots,A_{i-1},A_m,\ldots,A_i,   B' \cup A_2' $ is a red $\p_m$, here $A_{i-1}\cup A_m$ is red because the assumption of this lemma,
which is a contradiction. \hfill $\square$

The next lemma will tell us that the combination of two lemmas above forces a blue $\p_3$ under certain conditions.

\begin{lemma} \label{l:lm3}
If the edge $A_i \cup A_j$  is red  for all $1 \leq i \not = j
  \leq m$ and there is no red $\p_m$,  then there must be a blue $\p_3$.
\end{lemma}
\noindent
{\bf Proof:} We have two cases depending on the parity of $s$.

\begin{description}
\item [Case 1:] $s$ is even. We first show  that there is  at least one blue edge of type $(s/2,s/2,s)$. Suppose all edges of type $(s/2,s/2,s)$ are red. We pick two disjoint $s/2$-subsets $B'$ and $B''$ of $B$. Let $A_2'$ be an $s/2$-subset of $A_2$.  We define
\[
g_1=B' \cup A_2' \cup A_1, \quad  g_2=B''\cup \left(A_2 \setminus A_2'\right) \cup A_1,  \textrm{ and }  g_3=B'' \cup \left(A_2 \setminus A_2'\right) \cup  A_3.
\]
We get all $g_1, g_2,$ and $g_3$ are red since each of them is of type $(s/2,s/2,s)$. Now we get a red $\p_m$ with edges $g_1,g_2,g_3,f_3,\ldots, f_{m-1}$, which is a contradiction . Let $j$ be the smallest integer such that $1 \leq j \leq s/2$ and there is one blue edge of type $(j,s-j,s)$. It is easy to see $j$ is well-defined.  We observe that all edges of type $(s,0,s)$ are blue as we assume there is no red $\p_m$. If $j=1$, then a blue $\p_3$ is given by Lemma \ref{l:lm1} with $l=1$. If $j \geq 2$, then we get all edges of type $(j-1,s-j+1,s)$ are red by the minimality of $j$. Applying Lemma \ref{l:lm2} with $l=j-1$,  we get that all edges of type $(s-j+1,j-1,s)$ are blue. Now Lemma \ref{l:lm1} with $l=j$ gives us a blue $\p_3$.

\item[Case 2:]  $s$ is odd. We first consider the case where all edges of type $(j,s-j,s)$ are red for each $1 \leq j \leq \tfrac{s-1}{2}$. Using Lemma \ref{l:lm2} with $l=\tfrac{s-1}{2}$, we get all edges of type $(\tfrac{s+1}{2},\tfrac{s-1}{2},s)$ are blue. Let $B'$ be a $\tfrac{s+1}{2}$-subset of $B$, $A_2'$ and $A_2''$ be two disjoint $\tfrac{s-1}{2}$-subsets of $A_2$. We define
\[
g_1=B'\cup A_2'\cup A_1, \quad g_2=\left(B \setminus B'\right) \cup A_2'' \cup A_1,  \textrm{ and } g_3=\left(B \setminus B'\right) \cup A_2'' \cup A_m.
\]
 Because $g_1,g_2,$ and $g_3$ are of type  $(\tfrac{s+1}{2},\tfrac{s-1}{2},s)$, all of them are blue. We get a blue $\p_3$ with edges $g_1,g_2$, and $g_3$. If there is a blue edge of type $(j,s-j,s)$ with $1 \leq j \leq \tfrac{s-1}{2}$, then we repeat the argument in Case 1 to  get a blue $\p_3$.
 \end{description}
 \hfill $\square$

With all lemmas in hand, we are ready to prove Theorem \ref{t:1}.

\noindent
{\bf Proof of Theorem \ref{t:1}:} We will prove the theorem by induction on $n$. The base case is $n=3$.
Let $c$ be a red-blue coloring of edges of $K_{4s+1}^{2s}$. Since $4s+1 \geq R(\p_3, \p_2)$, either there is some red  $\p_3$, or there is some blue $\p_{2}$. We need only consider the latter case. We assume  a maximum blue path is $A_1, A_2, A_3$, where $|A_i|=s$ for each $1 \leq i \leq 3$.  Let $B$ be the remaining $s+1$ vertices and $B'$ be an arbitrary $s$-subset of $B$. Observe that  the edges $B' \cup A_1$ and $B' \cup A_3$ must be red as the maximum length of a blue path is two. If $A_1 \cup A_3$ is a blue edge, then a red $\p_3$ follows from Lemma \ref{l:lm3} by swapping colors. If $A_1 \cup A_3 $ is red, then $B'\cup A_1, A_1 \cup A_3, A_3 \cup B'$ form a red $\c_3$ for some $B' \subseteq B$. If there is no red $\p_3$, then there has to be a blue $\p_3$ by  Lemma \ref{l:lm3}, which is a contradiction. In either case, we are able to find a red $\p_3$ and we completed the proof for the base case.

Assume Theorem \ref{t:1} holds for all $3 \leq n \leq m-1 $ with  $m \geq 4$.
 Consider a red-blue coloring  of edges in $K_{(m+1)s+1}^{2s}$. Since $(m+1)s+1 \geq R(\p_{m-1},\p_3)=ms+1$ by  the inductive hypothesis, either  there is a red $\p_{m-1}$ or there is a blue $\p_3$. We need only consider the case in which the maximum length of a red path is $m-1$. Let $f_1,f_2,\ldots,f_{m-1}$ be a red $\p_{m-1}$, where $f_i=A_i \cup A_{i+1}$  for $1 \leq i \leq m-1$ and $|A_i|=s$ for each $1 \leq i \leq m$. Let  $B$ be the remaining $s+1$ vertices. Since there is no red $\p_m$, edges $B' \cup A_1$ and $B' \cup A_m$ must be blue for each $s$-subset $B'$ of $B$.
 We have the following mutually disjoint cases.
 \begin{description}
 \item[Case 1:]  Either  $A_1 \cup A_j$ is blue  for some $3 \leq j \leq m-1$ or $A_k \cup A_m$ is blue for some $2 \leq k \leq m-2$. Pick an arbitrary $s$-subset $B'$ of $B$.  We observe that $A_m, B',  A_1,  A_j$ form a blue $\p_3$ in the former case, and $A_1, B',  A_m,  A_k$ form a blue $\p_3$ in the latter case.
 \item[Case 2:]   $A_1 \cup A_i$  is red  for each $3 \leq i \leq m-1$ and $A_i \cup A_m$ is red for each $2 \leq i \leq m-2$. Moreover, there are $2 \leq j < k \leq m-1$ such that $k >j+1$ and $A_j \cup A_k$ is blue. We consider a new red $\p_{m-1}$ which is formed by
$A_1,A_2,\ldots,A_{k-1},A_m,A_{m-1},\ldots,A_k$.
Now $A_k$ is an ending $s$-set of this new path and we can find a blue $\p_3$ in the same way as in Case 1.
 \item[Case 3:]  We have  $A_i \cup A_j$ is red for all $1 \leq i \not = j \leq m$ such that $\{i,j\} \not = \{1,m\}$. Now if $A_1 \cup A_m$ is blue, then we can find a blue $\p_3$ by the same argument as Case 2.  Namely, we  find a new red $\p_{m-1}$ with one of $A_1$ and $A_m$ as an ending $s$-set but not the other one.   If $A_1 \cup A_m$ is red, then a blue $\p_3$ is ensured by Lemma \ref{l:lm3}.
 \end{description}
 \hfill $\square$

\section{Proof of Theorem \ref{t:2}}
For  a fixed $s \geq 1$, we will also prove Theorem \ref{t:2} by induction on $n$. Since the proof of the base case and the inductive step are similar, we sketch the idea for proving the inductive step here. We assume $R(\p_n,\p_4)=(n+1)s+1$ for all $4 \leq n \leq  m-1$.  For the inductive step,  let $c$ be a red-blue coloring of the edges of $K_{(m+1)s+1}^{2s}$. Since $(m+1)s+1 \geq R(\p_{m-1},\p_4)=ms+1$ by the inductive hypothesis,   either there is some red $\p_{m-1}$ or there is some blue $\p_4$. There is nothing to show if either there is some red $\p_m$ or a blue $\p_4$. Thus we assume that the maximum length of a red path is $m-1$; our goal is to find a blue $\p_4$ under this condition.  Let $A_1,A_2,\ldots,A_m$ be a fixed red $\p_{m-1}$ induced by $c$, where $\{A_1,A_2,\ldots,A_m\}$ is a collection of mutually disjoint $s$-sets of $[ms]$ and $f_i=A_i \cup A_{i+1}$ for $1 \leq i \leq m-1$. Let $B=[(m+1)s+1]\setminus [ms]$.   We will frequently replace some edges of the existing red $\p_{m-1}$ to obtain a new red $\p_{m-1}$ with new ending $s$-sets.
We will find that a blue edge $f$ with vertices from   $\cup_{i=1}^m A_i$ will help us to obtain a blue $\p_4$.
  There are many possible arrangements of the vertices of $f$. The simplest case is $f=A_i \cup A_j$ for some $1 \leq i \not =j \leq m$. We will show that we can always reduce the case where $f=A_i \cup A_j$ to the case where $f=A_1 \cup A_p$ for some $3 \leq p\leq m-1$.  If $f=A_1 \cup A_p$ is red, then the following lemmas tell us  how can we find the desired blue $\p_4$ under some conditions. We will repeatedly use the following fact:
\begin{fact}\label{f:f1}
Let $A_1,A_2,\ldots,A_m$ be a maximum red $\p_{m-1}$ induced by $c$. Then both $A_1 \cup B'$ and  $A_m \cup B'$ are blue for each $s$-subset $B'$ of $B$.
\end{fact}
The fact follows from the maximality of the red path $\p_{m-1}$.
%

For a fixed  red path $A_1,A_2,\ldots,A_m$, we say an edge $f$ is of type $(l,s-l,s)$ if $|f \cap  B|=l$, $|f \cap A_2|=s-l$, and   $A_p \subset f $ for some $1 \leq p \leq m$ with $p \not =2$.   Lemmas \ref{l:lm4}, \ref{l:lm5}, and \ref{l:lm6} play the same roles as   Lemmas \ref{l:lm1}, \ref{l:lm2}, and \ref{l:lm3}.
\begin{lemma} \label{l:lm4}
 Assume $A_1 \cup A_p$ is blue for some $3 \leq p \leq m-1$, $A_1 \cup A_{i}$ is red for all $ 3  \leq i \not = p \leq m-1$,  and $A_j \cup A_m$ is red for all $ 2 \leq j \leq m-2$.  Furthermore, assume there is no red $\p_m$.  For a fixed  $1 \leq l \leq \lfloor \tfrac{s}{2} \rfloor$,  if  all edges of type $(s+1-l,l-1,s)$ are blue and there is  a blue edge of type $(l,s-l,s)$,   then there exists  a  blue $\p_4$.
\end{lemma}
\noindent
{\bf Proof:} We assume that there is some blue edge $g_1$ of type $(l,s-l,s)$, say $g_1=B' \cup A_2' \cup A_j$,  where    $B'$ is an  $l$-subset of $B$,   $A_2'$ is an $(s-l)$-subset of $A_2$, and $j \not =2$.    We define $A_2''$ to be an arbitrary $(l-1)$-subset of $A_2 \setminus A_2'$.
We have two cases.
\begin{description}
\item[Case 1:] $j \in \{1,p\}$. Without loss of generality, we assume $j=p$. We define
\[
g_2=A_p \cup A_1, \quad g_3= \left( B\setminus B' \right) \cup  A_2'' \cup A_1, \textrm{ and } g_4= \left( B \setminus B' \right) \cup A_2'' \cup A_m.
\]
As $g_3$ and $g_4$ are of type $(s+1-l,l-1,s)$, both $g_3$ and $g_4$ are blue by assumption. The edge $g_2$ is also blue by the assumption. Now,  $g_1,g_2,g_3$, and $g_4$ form a blue $\p_4$. When $j=1$, we set $g_3= \left( B\setminus B' \right) \cup  A_2'' \cup A_p$. We still get a blue $\p_4$ with edges $g_1,g_2,g_3$, and $g_4$.
\item[Case 2:]  $j \not \in \{1,p\}$.   We define
\[
g_2=\left(B\setminus B'\right) \cup  A_2'' \cup A_j, \quad g_3=\left( B\setminus B' \right) \cup  A_2'' \cup A_1, \textrm{ and } g_4= A_1 \cup A_p.
\]
Both $g_2$ and $g_3$ are blue since they are of type $(s+1-l,l-1,s)$.  The assumption tells us $g_4$ is also blue. Now,
we obtain a blue $\p_4$ with edges $g_1,g_2,g_3$, and $g_4$.
\end{description}

 \hfill $\square$

We also have the following lemma which is similar to Lemma \ref{l:lm2}.
\begin{lemma} \label{l:lm5}
 Assume $A_1 \cup A_p$ is blue for some $3 \leq p \leq m-1$, $A_1 \cup A_{i}$ is red for all $ 3  \leq i \not = p \leq m-1$,  and $A_j \cup A_m$ is red for all $ 2 \leq j \leq m-2$.  Furthermore, assume there is no red $\p_m$.
For a fixed  $1 \leq l \leq \lfloor \tfrac{s}{2} \rfloor$, if  all edges of type $(l,s-l,s)$ are red, then  all edges of type $(s-l,l,s)$ are blue.
\end{lemma}
\noindent
{\bf Proof:} Suppose that there is an edge $g_1$ of type $(s-l,l,s)$ which is red.  We can assume $g_1=B' \cup A_2' \cup A_j$, where  $B'$ is a subset of $B$ with size $s-l$, $A_2'$ is a subset of $A_2$ with size $l$, and   $j \not =2$.
We first assume $j \not \in \{1,3\}$. Let $B''$ be an arbitrary $l$-subset of $B \setminus B'$. We get both $B'' \cup \left(A_2 \setminus A_2' \right) \cup A_1$ and $B'' \cup \left(A_2 \setminus A_2' \right) \cup A_3$ are red since both of them are of type $(l,s-l,s)$. Now,
$A_1,B'' \cup \left(A_2 \setminus A_2' \right), A_3,\cdots,A_{j-1},A_m,\ldots,A_j, B' \cup A_2'$ is a red $\p_m$, here $A_m \cup A_{j-1}$ is red  by assumption, which is   a contradiction.
For $j \in \{1,3\}$, we can find a red $\p_m$ similarly, see the proof of Lemma \ref{l:lm2}.  Therefore,  all edges of type $(s-l,l,s)$ must be blue.
 \hfill $\square$

The next lemma shows how can we get a blue $\p_4$ under some conditions.
\begin{lemma} \label{l:lm6}
 Assume $A_1 \cup A_p$ is blue for some $3 \leq p \leq m-1$, $A_1 \cup A_{i}$ is red for all $ 3  \leq i \not = p \leq m-1$,  and $A_j \cup A_m$ is red for all $ 2 \leq j \leq m-2$.  Furthermore, assume there is no red $\p_m$. Then there must be a blue $\p_4$.
\end{lemma}
\noindent
{\bf Proof:} Since the proof of this lemma uses the same idea as the one in the proof of Lemma \ref{l:lm3},  we outline it here.

 When $s$ is even,  we
define $j$ to be the smallest integer such that $1 \leq j \leq s/2$ and there is an edge of type $(j,s-j,s)$ which is blue. We are able to show $j$ is well-defined. If $j=1$ then a blue $\p_4$ is given by Lemma \ref{l:lm4} with $l=1$. If $j \geq 2$,  then we get all edges of type $(j-1,s-j+1,s)$ are red. Lemma \ref{l:lm5} with $l=j-1$ gives us that all edges of type $(s-j+1,j-1, s)$ are blue. Using Lemma \ref{l:lm4} with $l=j$, we can get a blue $\p_4$.

When $s$ is odd, we first show that there is a blue $\p_4$ if all edges of type $(j,s-j,s)$ are red for all $1 \leq j \leq \tfrac{s-1}{2}$.  Next we assume that there is a $1 \leq j \leq \tfrac{s-1}{2}$ such that  there is  a blue  edge   of type $(j,s,s-j)$. We can get a blue $\p_4$ using Lemma \ref{l:lm4} and Lemma \ref{l:lm5} as we did  for proving Lemma \ref{l:lm3}.
 \hfill $\square$

We have the following lemma for the special case where $m=5$.

\begin{lemma}\label{m:m1}
Assume $A_1,A_2,\ldots,A_5$ is a red $\p_4$ and there is no red $\p_5$. If $A_1 \cup A_3$  and $A_3 \cup A_5$ are blue,    then there is a blue $\p_4$.
\end{lemma}

\noindent
{\bf Proof:}
 We first show that we need  only  consider the case in which both $A_1 \cup A_4$ and $A_2 \cup A_5$ are red.
Let $B'$ be an $s$-subset of $B$ and recall Fact  \ref{f:f1}. We get a blue $\p_4=A_4,A_1,A_3,A_5,B'$  if  $A_1 \cup A_4$  is  blue. Similarly,   $A_2,A_5,A_3,A_1,B'$ is  a blue $\p_4$ if $A_2 \cup A_5$ is blue.
We have two cases depending on the color of $A_1 \cup A_5$.

\begin{description}
 \item [Case 1:] The edge $A_1 \cup A_5$ is red.   We form a new red $\p_4$ as $A_1,A_5,A_2,A_3,A_4$. We note that $B' \cup A_4$ is blue as Fact \ref{f:f1} and $A_5,A_3,A_1,B',A_4$ form a blue $\p_4$.

 \item [Case 2:] The edge $A_1 \cup A_5$ is blue. For an edge $f$ of type $(l,s-l,s)$, we define the {\em center} of $f$ to be the unique $A_i$ such that $A_i \subset f$.  Fact 1 implies that all edges of type $(s,0,s)$ with center from $\{A_1,A_3,A_5\}$ must be blue.
 We have the following two claims which are similar to Lemma \ref{l:lm4} and Lemma \ref{l:lm5}.
 \end{description}

 \noindent
 {\bf Claim 1:} For a fixed $1 \leq l \leq \lfloor \tfrac{s}{2} \rfloor$, if all edges of type $(s+1-l,l-1,s)$ with center from $\{A_1,A_3,A_5\}$ are blue and there is  a blue edge of type $(l,s-l,s)$  with center from $\{A_1,A_3,A_5\}$, then there exists a blue $\p_4$.

{\em Proof of Claim 1:} By  symmetry of $A_1$, $A_3$, and $A_5$,  we can assume $B' \cup A_2' \cup  A_1$ is blue, where $B'$ is an $l$-subset of $B$ and $A_2'$ is an $(s-l)$-subset of $A_2$.   Let $A_2''$ be an $(l-1)$-subset of $A_2 \setminus A_2'$. We get both  $(B\setminus B') \cup A_2'' \cup  A_1$ and $(B\setminus B') \cup A_2'' \cup  A_3$ are blue as they are of type $(s+1-l,l-1,s)$.  We note that  $B' \cup A_2', A_1, (B\setminus B') \cup A_2'', A_3,A_5$ form a blue $\p_4$.

\noindent
{\bf Claim 2:} For each fixed $1 \leq l \leq \lfloor \tfrac{s}{2} \rfloor$,  if all edges of type $(l,s-l,s)$ with center from $\{A_1,A_3,A_5\}$ are red,  then  all edges of type  $(s-l,,l,s)$ with center from $\{A_1,A_3,A_5\}$ must be blue.

{\em Proof the Claim 2:} By  symmetry of $A_1$, $A_3$, and $A_5$, we can  assume $g_1=B' \cup A_2' \cup A_1$ is red, where $B'$ is an $(s-l)$-subset of $B$ and $A_2'$ is an $l$-subset of $A_2$.
 Pick an $l$-subset $B''$ of $B\setminus B'$. We define
\[
g_2=B'' \cup (A_2\setminus A_2') \cup A_1 \textrm{ and } g_3=B'' \cup  (A_2 \setminus A_2') \cup A_3.
\]
We get $g_2$ and $g_3$ are of type $(l,s-l,s)$ and they are red by assumption. Now $B'\cup A_2', A_1, B'' \cup (A_2 \setminus A_2'), A_3, A_4, A_5$ form a red $\p_5$, which is a contradiction.

To find a blue $\p_4$, we repeat  the argument in the proof of Lemma \ref{l:lm6} as follows depending on the parity of $s$. If $s$ is even, then we first show there is one edge of type $(s/2,s/2,s)$ with center from $\{A_1,A_3,A_5\}$ which is blue. Suppose not. Let $B'$ and $B''$ be two disjoint $s/2$-subsets of $B$ and $A_2'$ be an $s/2$-subset of $A_2$. Now $B'\cup A_2', A_1, B'' \cup (A_2 \setminus A_2'), A_3,A_4,A_5$ form a red $\p_5$ and we get a contradiction.
We define $j$ as the smallest  integer such that $1 \leq j \leq s/2$  and  there is a blue edge of type $(j,s-j,s)$ with center from $\{A_1,A_3,A_5\}$. Then $j$ is well-defined. If $j=1$, then   all edges of type $(s,0,s)$ with center from $\{A_1,A_3,A_5\}$ are blue. A desired  blue $\p_4$ is given by Claim 1.  If $j \geq 2$, then
we get the assumption in Claim 2 for $l=j-1$. The conclusion of Claim 2 with $l=j-1$ together with the definition of $j$ give us the assumption in Claim 1 with $l=j$. Now Claim 1 with $l=j$ gives us a required blue $\p_4$.

If $s$ is odd, then defining $j$ similarly, we can show $j \leq \tfrac{s-1}{2}$.
Repeating the argument for the case where $s$ is even,  we can get a blue $\p_4$.  \hfill $\square$

As we mentioned before, the existence of a blue edge $f=A_i \cup A_j$ is helpful for finding a blue $\p_4$.  The next lemma will show the case in which $f=A_1 \cup A_p$ for some $3 \leq p \leq m-1$.
\begin{lemma}\label{l:lm7}
If $A_1 \cup A_p$ is blue for some $3 \leq p \leq m-1$, then there is a blue $\p_4$.
\end{lemma}
{\bf Proof:} If  there is some   $2 \leq j \not =p \leq m-2$ such that  $A_j \cup A_m$ is  blue, then we take   an $s$-subset $B'$ of $B$.
 Fact \ref{f:f1} implies that both $A_1 \cup B'$ and $B' \cup A_m$ are blue. Note that $A_p,A_1,B',A_m,A_j$ form a blue $\p_4$.  In the remaining proof, we  assume $A_j \cup A_m$ is red for each $2 \leq j \not = p  \leq m-2$.
Note that the above argument gives us the assumptions in Lemma \ref{l:lm6} for $m=4$; thus a desired blue $\p_4$ is ensured by Lemma  \ref{l:lm6} for $m=4$.

 We first consider the case where $m \geq 6$. We get that either $p-1 \geq 3$ or $m-p \geq 3$.  We wish to show that it suffices to consider the case where $A_p \cup A_m$ is red. Suppose $A_p \cup A_m$ is blue. We aim to find a blue $\p_4$ directly.  The  idea is that we find a new red path with length $m-1$ which contains $A_q$ as an ending $s$-set for some $q \not \in \{1,p,m\}$.
If $p-1 \geq 3$, then we consider a new  red path $A_1,A_2,A_m,A_3,\ldots,A_p, \ldots,A_{m-1}$.  If $m-p \geq 3$, then we look at  a new red path $A_1,A_2,\ldots,A_p,\ldots,A_{m-2},A_{m},A_{m-1}$.
Here we use the assumption $A_j \cup A_m$ is red for each $2 \leq j \not = p \leq m-2$.
Fact \ref{f:f1} implies $A_1 \cup B'$ and  $A_{m-1} \cup B'$ are blue for each $s$-subset $B'$ of $B$. Now, $A_{m},A_p,A_1,B',A_{m-1}$ is a blue $\p_4$ in both cases.
Thus, we can assume $A_p \cup A_m$ is red.

Under the assumption $A_1 \cup A_p$ is blue and $A_j \cup A_m$ is red for each $2 \leq j \leq m-2$,  we wish to show that it is sufficient to examine the case where $A_1 \cup A_j$  is red for each $3 \leq j \not = p \leq m-1$.  Suppose $A_1 \cup A_j$  is blue for some $3 \leq j \not = p \leq m-1$.  If $j<p$, then we consider a new red path $A_1,\ldots,A_j,\ldots,A_{p-1},A_m,\ldots, A_p$. If $j>p$, then we form a new red $\p_{m-1}$ as $A_1,\ldots, A_p, \ldots, A_{j-1}, A_m, \ldots, A_j$.
Take an $s$-subset $B'$ of $B$. Then
$A_m, B', A_p, A_1, A_j$ will be a blue $\p_4$ in the first case and  $A_m, B', A_j, A_1, A_p$ is a blue $\p_4$  in the second case.
We get the assumptions stated in Lemma \ref{l:lm6} and   a desired blue $\p_4$ is given by Lemma \ref{l:lm6} for $m \geq 6$.

  Lastly, we prove the result for $m=5$. If $p=4$, then we first show that we can assume  $A_2 \cup A_5$ and $A_3 \cup A_5$ are red. Based on this assumption, we can  assume further $A_1 \cup A_3$ is red. Since the argument here is exactly the same as the case $m \geq 6$,  it is omitted.  If $p=3$ and $A_3 \cup A_5$ is red, then we can show that we need only consider the case where $A_2 \cup A_5$ and $A_1 \cup A_4$ are red by the same argument as the one for $m \geq 6$. We get the assumptions in Lemma \ref{l:lm6} for these two cases.
A blue $\p_4$ is given by Lemma \ref{l:lm6}.  If $p=3$ and $A_3 \cup A_5$ is blue,  then a blue $\p_4$ is given by Lemma \ref{m:m1}.
 \hfill $\square$

The next lemma will tell us  that we can reduce the general case where $f=A_i \cup A_j$ to the case where $f=A_1 \cup A_p$.

\begin{lemma} \label{l:lm8}
If there is some blue edge $f=A_i \cup A_j$ for some $1 \leq i \not = j \leq m$, then there is a blue $\p_4$.
\end{lemma}
\noindent
{\bf Proof:} We have the following mutually disjoint cases.
\begin{description}
\item[Case 1:] $|\{i,j\} \cap \{1,m\}|=1$. Note that the case where $f=A_j \cup A_m$ is the same as the case where $f=A_1 \cup A_p$ by  symmetry, so this case is proved by Lemma \ref{l:lm7}.
\item[Case 2:] $2 \leq i < j \leq m$ and $A_p \cup A_q$ is red for all $|\{p,q\} \cap \{1,m\}|=1$. We observe that $A_j, \ldots,A_m,A_{j-1}, \ldots,A_i, \ldots, A_1$ is a new red $\p_{m-1}$ and we can reduce it to Case 1.
\item[Case 3:] $\{i,j\}=\{1,m\}$ and $A_p \cup A_q$ is red  for all $\{p,q\} \not = \{1,m\}$. We form a new red $\p_{m-1}$ as $A_1,A_2,A_m, \ldots ,A_3$  and we reduce it to Case 1.
\end{description}
\hfill $\square$

We already showed how a blue edge $f=A_i \cup A_j$ helped us to get a blue $\p_4$. In general, $f$ could intersect more than two $A_i$'s. The next lemma will give us a  blue $\p_4  $ for other possible intersections between $f$ and $A_i$'s. We first introduce some related notation.
   Given    a red path  $\p_{m-1}=A_1,A_2, \ldots, A_m$ and an edge $f$ with $f \subseteq \cup_{i=1}^m A_i$, let $S(\p_{m-1},f)=\{i: 1 \leq i \leq m \textrm{ and } f \cap A_i \not = \emptyset\}$. We say a fixed coloring $c$ has \pr(i) if the existence of  some edges $f$ and a red path $\p_{m-1}$ satisfying $S(\p_{m-1},f)=i$ implies the existence of  a blue $\p_4$.
We have the following lemma.
\begin{lemma} \label{l:lm9}
For a fixed red-blue coloring $c$ of edges of $K_{(m+1)s+1}^{2s}$ without a red $\p_m$, then the coloring $c$ has \pr(i) for each $2 \leq i \leq \min\{m,s\}$.
\end{lemma}
\noindent
{\bf Proof:} We proceed by induction on $i$. The base case where $i=2$ is given by Lemma \ref{l:lm8}.  We  assume $c$ has \pr(i) for all $2 \leq i \leq k-1$. For the inductive step, let us fix a red $\p_{m-1}$ and a blue edge $f$ satisfying $|S(\p_{m-1},f)|=k$.
 Without loss of generality, we assume $S(\p_{m-1},f)=\{1,\ldots,k\}$. Let
 $A_i'=f\cap A_i$ for each $1 \leq i \leq k$.

If  $k \geq 4$ then we can assume $|A_1'| \leq |A_2'| \leq \cdots \leq |A_k'|$. Clearly, $| A_1' \cup A_2'| \leq s $ by the pigeonhole principle. Let $C$ be a subset of $A_1 \cup A_2$ such that $A_1' \cup A_2' \subseteq C$ and $|C|=s$.
If the edge $A_3 \cup C$ is blue, then a  blue $\p_4$ is give by the inductive hypothesis by noticing $|S(\p_{m-1}, A_3 \cup C)|=3$.
Thus we can assume $A_3 \cup C$ is red.
Let $C'=(A_1 \cup A_2) \setminus C$ and we consider a new red path $\p_{m-1}'=A_m, \ldots, A_3,C,C'$, here $C \cup C'=A_1 \cup A_2$.
 We get  a blue $\p_4$ by the inductive hypothesis as $|S(\p_{m-1}',f )|=k-1$.

If $k=3$, then additional arguments are needed.
 We can  assume $A_i \cup A_j$ is red for all $1 \leq i \not =j \leq m$; otherwise the base case gives us a blue $\p_4$.

 We first consider that there is some $1 \leq i \leq 3 $ such that $A_i \subseteq f$. Without loss of generality, we assume $A_1 \subseteq f$.  Let
 $C=A_2' \cup A_3'$ and $C'=(A_2 \cup A_3) \setminus C$. We define
 $g_2=C' \cup  A_m$. If $g_2$ is blue, then let $g_3=A_m \cup B'$ and $g_4=B' \cup A_1$, here $B' \subseteq B$ and $|B'|=s$.   Fact \ref{f:f1} implies $g_3$ and $g_4$ are blue. Now, $g_2,g_3,g_4,f$ form a blue $\p_4$. If $g_2$ is red, then we form a new red path $\p'_{m-1}=C,C',A_m,A_{m-1},\ldots,A_4,A_1 $. Note $|f \cap \p'_{m-1}|=2$ and the base case gives us a blue $\p_4$.

  If $|A_i \cap f| < s$ for each $1 \leq i \leq 3$, then  we observe  $|A_i' \cup A_j'| \geq  s$ for some $1 \leq i \not =j \leq 3$ by the pigeonhole principle.  We assume $|A_2' \cup A_3'| \geq s$ and  pick a subset $A_2''$ of $A_2'$ such that $|A_2'' \cup A_3'|=s$. Let $C=A_2'' \cup A_3'$
 and    $C'=(A_2 \cup A_3)  \setminus C$.
   We need only consider the case where $C' \cup A_4$ and  $A_1 \cup C$ are red. If $g=C' \cup A_4$ is blue, then a blue $\p_4$ is given by the  previous case by observing $|g \cap \p_{m-1}|=3$ and $A_4 \subset g$.  We have a similar argument for $A_1 \cup C$.  When both $C' \cup A_4$ and $A_1 \cup C$ are red,
   we observe $C',A_4,\ldots,A_m,A_1,C$ is a new red $\p_{m-1}'$, $|S(\p'_{m-1},f)|=3$, and $C \subseteq f$.  We reduce this case to the previous case. \hfill $\square$

We already know how to find a blue $\p_4$ if there is some blue $f$ such that $f \subseteq \cup_{i=1}^m A_i$.
Next, we assume $f$ is red for all $f \subseteq \cup_{i=1}^m A_i$ and show how can we   find a blue $\p_4$ under this assumption. We need one more definition. Fix a red path $A_1,A_2,\cdots, A_m$ and let $B$ be the remaining $s+1$ vertices. For each $1 \leq l \leq s$,   we say $f$ is of type $(s-l,s+l)$ if $|f \cap B|=s-l$ and $|f\cap (\cup_{i=1}^m A_m)|=s+l$.

\begin{lemma}\label{l:lm10}
Let $A_1,A_2,\ldots, A_m$ be a red $\p_{m-1}$.  Assume all edges $f$ with $f \subseteq \cup_{i=1}^m A_i$ are red and there is no red $\p_m$. For each $1 \leq l \leq \lfloor \tfrac{s}{2} \rfloor$, if all edges of type $(s-l+1,s+l-1)$ are blue and there is a red edge of type $(s-l,s+l)$, then there exists    a blue $\p_4$.
\end{lemma}
\noindent
{\bf Proof:} Suppose that  there is a red edge $f$ of type $(s-l,s+l)$. Without loss of generality, we can assume $f=B' \cup A_1 \cup A_2'$,
here $A_2' \subseteq A_2$ with $|A_2'|=l$ and $B' \subseteq B$ with $|B'|=s-l$.  Let $B''$ be an $l$-subset of $B \setminus B'$. We define
\[
g_1=B'' \cup (A_2 \setminus A_2')\cup A_3 \textrm{ and } g_2=B'' \cup ( A_2 \setminus A_2') \cup  A_4.
\]
We get both $g_1$ and $g_2$ are blue. Otherwise, if $g_1$ is red, then $B''\cup (A_2\setminus A_2' ), A_3,  \ldots, A_m, A_1, B'\cup A_2'$ is a red $\p_m$, which is a contradiction.  If $g_2$ is red, then we can find a contradiction similarly.  Let $A_2''$ be an $(l-1)$-subset of $A_2'$. We  define
\[
g_3=(B \setminus B'')  \cup A_2'' \cup A_1 \textrm{ and } g_4=A_2'' \cup (B \setminus B'') \cup  A_3.
\]
We observe that both $g_3$ and $g_4$ are of type $(s-l+1,s+l-1)$. Thus both $g_3$ and $g_4$ are blue by assumption. Now, $g_3,g_4,g_1,g_2$ is a blue $\p_4$.  \hfill $\square$

The next lemma will show how  Lemma \ref{l:lm10} guarantees a blue $\p_4$.
\begin{lemma}\label{l:lm11}
Let $A_1,A_2,\ldots, A_m$ be a red $\p_{m-1}$.  Assume all edges $f$  satisfying $f \subseteq \cup_{i=1}^m A_i$ are red and there is no red $\p_m$. Then we have a blue $\p_4$.
\end{lemma}
\noindent
{\bf Proof:}  Let $j$ be the smallest integer such that $1 \leq j \leq \lfloor \tfrac{s}{2} \rfloor$ and there is a red edge of type $(s-j,s+j)$. If there is such a $j$ then we get all edges of type $(s-j+1,s+j-1)$ are blue by the choice of $j$. In the case where $j=1$,   all edges of type $(s,s)$ are blue by Fact 1.  Applying Lemma \ref{l:lm10} with $l=j$, we get a blue $\p_4$.  If there is no such a $j$ then all edges of type $\left(s-\lfloor \tfrac{s}{2} \rfloor,s+\lfloor \tfrac{s}{2} \rfloor\right)$ are blue.

When $s$ is odd, let $B'$ be a  subset of $B$ with size $\tfrac{s+1}{2}$.  Let  $A_1'$ and  $A_1''$ be two disjoint subsets of  $A_1$ with  size $\tfrac{s-1}{2}$. We define
\[
g_1=A_2 \cup A_1' \cup B', \quad g_2=A_1' \cup B' \cup A_3, \quad g_3=A_3 \cup A_1'' \cup  (B\setminus B'), \textrm{ and } g_4=A_1''\cup (B \setminus B') \cup A_4.
\]
We observe  $g_1,g_2,g_3,g_4$ form a blue $\p_4$ as each $g_i$ is of type $(\tfrac{s+1}{2},\tfrac{3s-1}{2})$ for each $1 \leq i \leq 4$.

When $s$ is even,  let $B'$ and $B''$ be two disjoint  subsets $B$ with size $\tfrac{s}{2}$  and  $A_2'$  be a  subset of  $A_2$ with  size $\tfrac{s}{2}$. We define
\[
g_1=A_1 \cup A_2' \cup B', \quad g_2=A_2' \cup B' \cup A_3, \quad g_3=A_3 \cup (A_2 \setminus A_2') \cup  B'', \textrm{ and } g_4=(A_2 \setminus A_2') \cup  B'' \cup A_4.
\]
We notice that  $g_1,g_2,g_3,g_4$ form a blue $\p_4$ as each $g_i$ is of type $(s/2,3s/2)$ for $1 \leq i \leq 4$. In either case, we are able to find a blue $\p_4$.
\hfill $\square$

We are now ready to prove Theorem \ref{t:2}.

\noindent
{\bf Proof of Theorem \ref{t:2}:}  We prove the theorem by induction on $n$. For the base case, let $c$ be a  red-blue coloring of edges in $K_{5s+1}^{2s}$.
As $5s+1 \geq R(\p_4,\p_3)$  by Theorem \ref{t:1}, either we have a red $\p_4$ or we have a blue $\p_3$.
 There is  nothing to show for the former case.  Thus we assume $A_1,A_2,A_3,A_4$ is a  maximum blue path. If there is a red edge with vertices from $\cup_{i=1}^4 A_i$, then we have a red $\p_4$ by Lemma \ref{l:lm9} with colors swapped. Otherwise, we switch colors in Lemma \ref{l:lm11} to get a red $\p_4$.

The inductive step is given by Lemma \ref{l:lm9} and Lemma \ref{l:lm11}.
 \hfill $\square$
\section{Concluding remarks}
In this paper, we give a partial affirmative answer to Question \ref{q:q1} for $s=r/2$, $r$  even, and $m \in \{3,4\}$. However, unlike in \cite{grs}, we are not able to determine the Ramsey number of small $r/2$-cycles for even $r$. A possible reason is following. The authors in \cite{grs} proved the following statement.  Let $c$ be a red-blue coloring of edges in $K^r_N$, here $N=(r-1)n+\lfloor \tfrac{m+1}{2}\rfloor$. If   $\c^{r,1}_n \subseteq \fr$, then either $\p^{r,1}_n \subseteq \fr$ or $\p^{r,1}_m \subseteq \fb$. Also,  if   $\c^{r,1}_n \subseteq \fr$, then either $\p^{r,1}_n \subseteq \fr$ or $\c^{r,1}_m \subseteq \fb$. The statement above is a very important fact for $s=1$; it helps to determine the values of $R(\p^{r,1}_n,\p^{r,1}_m)$, $R(\p^{r,1}_n,\c^{r,1}_m)$, and $R(\c^{r,1}_n,\c^{r,1}_m)$.
We can not prove a similar lemma for $s=r/2$ and $r$ even since after we fix a red $\c^{r,r/2}_n$,  no vertices remain.  It would be helpful to prove a lemma which connects $R(\p^{r,r/2}_n,\p^{r,r/2}_m)$ to $R(\c^{r,r/2}_n,\c^{r,r/2}_m)$.

To answer Question \ref{q:q1}, we need to determine the exact value of the Ramsey number of each type of path; it is very possible that we need different techniques to deal with different types of paths.
There are many other interesting questions on Ramsey number of paths and cycles in hypergraphs.   The only known results  addressing   tight cycles is due to Haxell et  al. \cite{h2} who examined the asymptotic value of $R(\c_n^{3,2},\c_n^{3,2})$.  A natural question is to determine the exact value of the Ramsey number of tight paths and  cycles.

\vspace{0.3cm}

\noindent
{\bf Acknowledgement:} The author thanks  G.~R.~Omidi and anonymous referees  for their valuable comments on earlier versions of this paper. Part of
this work was done while the author was a postdoc at University of California, San Diego. The author thanks grants ONR MURI N000140810747 and AFSOR AF/SUB 552082.

\end{document}